\documentclass{ifacconf}
\usepackage{url}
\usepackage{amsmath}
\usepackage{graphicx}      
\usepackage{natbib}        
\usepackage{graphics} 
\usepackage{epsfig} 
\usepackage{mathptmx} 
\usepackage{times} 
\usepackage{amsmath} 
\usepackage{amssymb}  
\usepackage{multirow}
\input{epsf}
\usepackage{graphicx}
\usepackage{graphicx, subfigure}
\usepackage[utf8]{inputenc}
\usepackage[english]{babel}

\def\figurename{FIG}

\def\figurename{FIG}
\def\beq{\begin{equation}}
\def\eeq{\end{equation}}
\def\baq{\begin{eqnarray}}
\def\eaq{\end{eqnarray}}
\def\bal{\begin{align} }
\def\eal{\end{align} }
\def\bc{\begin{center}}
\def\ec{\end{center}}
\def\ds{\displaystyle}
\def\fai{\varphi}
\def\eps{\varepsilon}

\def\vtheta{\vartheta}

\def\fai{\varphi}
\def\eps{\varepsilon}

\def\vtheta{\vartheta}
\def\gr{\gamma_{r}}
\def\gr1{\gamma_{r1}}

\def\fai{\varphi}
\def\eps{\varepsilon}

\def\vtheta{\vartheta}

\def\figurename{FIG}

\def\figurename{FIG}
\def\beq{\begin{equation}}
\def\eeq{\end{equation}}
\def\baq{\begin{eqnarray}}
\def\eaq{\end{eqnarray}}
\def\bal{\begin{align} }
\def\eal{\end{align} }
\def\bc{\begin{center}}
\def\ec{\end{center}}
\def\ds{\displaystyle}
\def\fai{\varphi}
\def\eps{\varepsilon}

\def\vtheta{\vartheta}
\def\gr{\gamma_{r}}
\def\gr1{\gamma_{r1}}

\def\tt{\theta}
\def\fai{\varphi}
\def\eps{\varepsilon}

\def\vtheta{\vartheta}
\def\gr{\gamma_{r}}
\def\gr1{\gamma_{r1}}

\def\fai{\varphi}
\def\eps{\varepsilon}

\def\vtheta{\vartheta}

\begin{document}
\begin{frontmatter}

\title{Richardson Approach or Direct Methods ?  What to Apply in the Ill-Conditioned Least Squares Problem\thanksref{footnoteinfo}}

\thanks[footnoteinfo]{This work was not supported by any organization}

\author[First]{Alexander Stotsky}

\address[First]{Department of Computer Science and Engineering,
 Chalmers University and University of Gothenburg, SE-412 96 Gothenburg, Sweden
 \\ e-mail:  alexander.stotsky@chalmers.se }

\begin{abstract}                
: This report shows on real data that the direct methods such as LDL decomposition and Gaussian elimination for solving linear systems with ill-conditioned matrices provide inaccurate results due to divisions
by very  small numbers, which in turn results in peaking phenomena and large estimation errors.
Richardson iteration  provides accurate results without peaking phenomena since division by small numbers is
absent in the Richardson approach.
\\ In addition, two preconditioners are considered and compared in the Richardson iteration: 1) the simplest and robust preconditioner based on the maximum row sum matrix norm and  2) the optimal one based
on calculation of the eigenvalues. It is shown that the simplest preconditioner is more
robust for ill-conditioned case and therefore it is  recommended for many applications.
\end{abstract}


\end{frontmatter}
\section{Introduction}
\noindent
The solution of the system of linear equations  with {\it SPD (Symmetric and Positive Definite)} ill-conditioned matrix is required in many application areas such as control, system identification, signal processing,
statistics as well as in  many  big data applications.
\\ The matrix $A$ is SPD ill-conditioned (a) information matrix for systems with harmonic regressor
and short window sizes, \cite{cdc54}, \cite{sto2019}, (b) Gram matrix due to the squaring of the condition number in the least squares method, \cite{bjor}, \cite{ljung}, (c) mass matrix in finite element method, \cite{stic} and mass matrix (lumped mass matrix) for  mechanical systems with singular perturbations, \cite{mostafa},
(d) state matrix for systems of linear equations, \cite{has} and in many other applications.
Moreover, many mechanical and electrical systems are singulary perturbed (have modes with different time-scales), \cite{kok} and considered as stiff and ill-conditioned systems, which potentially extends application areas of  this model.
\\ Ill-conditioning implies robustness problems (sensitivity to numerical calculations) and
imposes additional requirements on the accuracy of the solution of algebraic equations.
The accuracy requirements is the motivation for application of the Richardson iteration which is driven by
the residual error and has filtering (averaging) property.  The residual error is smoothed and remains bounded for a sufficiently large step number in this iteration, where the bound
depends on inaccuracies, providing the best possible solution in finite digit calculations.
\\ The performance of Richardson algorithms strongly depends on the preconditioning and therefore development
of new computationally efficient preconditioners is required. To this end the movement of the window is presented as rank two update of the information matrix.  Indeed, the least squares estimation of the frequency contents  of oscillating signal
in the window of the size $w$ which is moving in time can be presented in the
following form:
\begin{align}
 A_k \theta_k &= b_k, ~~ b_k = \sum_{j = k -(w-1)}^{j=k} \fai_{j}~ y_j = b_{k-1} + d_k
 \label{lineq1} \\
  b_{k-1}  &= \sum_{j = k - w}^{j=k-1} \fai_{j}~ y_{j}, ~~ d_k =   \fai_{k} ~ y_k -   \fai_{k-w} ~ y_{k-w}
  \label{lineq2} \\
   A_k &= \sum_{j = k -(w-1)}^{j=k} \fai_{j}~ \fai_{j}^{T}  =  A_{k-1} + R_k \label{ak1} \\
   A_{k-1} &= \sum_{j = k - w}^{j=k-1} \fai_{j}~ \fai_{j}^{T}, R_k =  \fai_{k}~ \fai_{k}^{T}  -  \fai_{k-w}~ \fai_{k-w}^{T}   \label{rk1} \\
    \fai_k^{T} &=  [ cos( q_0 k ) ~ sin(q_0 k ) ~ ...~ cos( q_h k ) ~ sin( q_h k ) ]  \label{fai1}
\end{align}
where the oscillating signal ${y}_k$ is approximated using the model $\hat{y}_k = \fai_k^{T} \theta_k $
with the harmonic regressor  (\ref{fai1}), where $q_0,...q_h$ are the frequencies.
The parameter vector $\theta_k$ should be calculated in each step with desired accuracy
as the solution of the algebraic equation (\ref{lineq1}), which is associated with minimization of the following error
$\sum_{j = k -(w-1)}^{j=k} (y_j - \hat{y}_j)^2$, where ${y}_k = \fai_k^{T} \theta_* + \xi_k $ and  $\theta_* $ is
the vector of unknown parameters and $ \xi_k $ is the noise.
\\
The information matrix $A_k$ is defined in (\ref{ak1})
as the sum of rank one matrices and as the rank two update, $R_k$ of the matrix  $A_{k-1}$,
$k \ge w + 1$. Rank two update is associated with the movement of the window, where
the new data $\fai_{k}$, $y_k$ enter the window and the data $\fai_{k-w}$, $y_{k-w}$ leave the window
in step $k$.
\\ Ill-conditioning of the matrix $A_k$ implies robustness problems (sensitivity to numerical calculations) and
imposes additional requirements on the accuracy of the solution of (\ref{lineq1})
(which are especially pronounced in finite-digit calculations) since small changes in $b_k$ due to measurement, truncation, accumulation, rounding and other errors result in significant changes in $\theta_k$.
In addition, ill-conditioning implies slow convergence of the iterative procedures.
\begin{figure}
\centerline{\psfig{figure=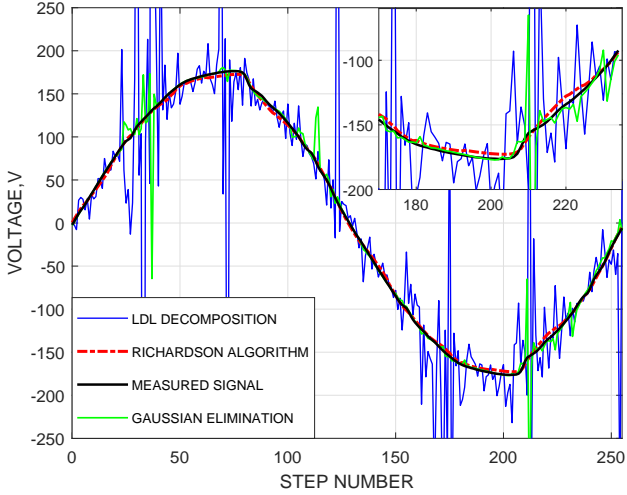,height=70mm}}
\begin{center}
\caption{\small{Measured voltage wave form  of a single cycle is
plotted with the black line.The Figure shows approximation performance of different types of
computational algorithms for estimation of the parameter vector.
Approximation with standard matrix inversion algorithm
(matrix inversion using LDL decomposition,implemented as standard routine in Matlab)
 is plotted with the blue line, and the approximation with the standard routine which realises
 Gaussian elimination method is plotted with the green line. Approximation with Richardson parameter estimation algorithm
 is plotted with the red dashed line.
}}
\label{figsurface3}
\end{center}
\end{figure}
\section{Which Algorithms? \\ Richardson or Direct Methods}
\noindent
{\it Richardson Iterations.} Richardson framework provides simple assessment and quantification of the
trade-off between the accuracy and computational burden, associated with the
concept of approximate computing and can be chosen as the
most promising solution for ill-conditioned systems.
The Richardson algorithms can be written in the following form:
\begin{align}
 \vtheta_i &= \vtheta_{i-1} - G_i ~ \{A_k \vtheta_{i-1} - b_k \},~F_i = I - G_{i} A_k,~ A_k \theta_k &= b_k  \label{ttk1}
\end{align}
where $\vtheta_i$ is the vector that estimates unknown parameters, $\theta_k = \vtheta_i $,
$G_i$ is associated with iterative matrix  inversion algorithms, which minimize inversion error $F_i$, where
$I$ is the identity matrix.
\\ The norm of the residual error
$\| A \tt_{k} - b \| \le \delta $, where $\delta >0$ is
a given bound can be used for a proper choice of the number of iterative steps providing pre-specified upper bound of the accuracy according to the concept of approximate computing.
\\ {\it Direct Methods.}
Gaussian elimination method is much faster and more accurate than the methods associated with the matrix inversion.
Gaussian method may produce residuals of the order of machine accuracy, but the solution is often not reliable
due to numerical instability, \cite{matlab}, \cite{dru}.
Namely, the solution is often accompanied with the peaking phenomena in the ill-conditioned case due to division by very small numbers, which is not present in the Richardson
approach, where the norm of the residual error can be kept uniformly within the  bound $\delta$.
\\ {\it Comparisons on Real Data.}
The one-phase synchronized voltage waveform measured at the wall outlet (approximately 120V RMS) is used for comparisons of Richardson and direct methods, \cite{ieee}. The sampling measurement rate is $256 $ points per cycle.Measured voltage wave form  of a single cycle is
plotted in Figure~\ref{figsurface3} with the black line.
The signal was approximated with the system
with harmonic regressor, (\ref{fai1}) which contains the fundamental frequency of $q_0=60$ Hz and four
higher harmonics. The wave form is approximated in the moving window of a relatively small size.
The condition number of the information matrix varies  significantly as the function
of step number and has the average order of $10^7$, indicating ill-conditioning.
Extreme ill-conditioning was detected
in several steps where the condition number reaches the order of $10^9$.
The Figure shows approximation performance of different types of
computational algorithms for estimation of the parameter vector.
Approximation with standard matrix inversion algorithm
(matrix inversion using LDL decomposition,implemented as standard routine in Matlab)
 is plotted with the blue line, and the approximation with the standard routine which realises
 Gaussian elimination method is plotted with the green line. Approximation with Richardson parameter estimation algorithm
 is plotted with the red dashed line.
 \\ Average accuracy (residual) of the method which is based on the matrix inversion
 is much worse than the average accuracies of Richardson and Gaussian method. Moreover, the matrix inversion method shows essential deterioration of the accuracy (due to division by very small numbers)
 in the large number points. Gaussian method provides very high average
 accuracy compared to the matrix inversion and Richardson methods.
 However, the accuracy deteriorates in some points and becomes worse than the accuracy of the Richardson method.
 The number of steps of the Richardson algorithm was optimized in each step of the moving window
to guarantee the pre-specified upper bound of the accuracy, which eliminated peaking phenomena
and reduced computational time.
Average accuracy of the Richardson algorithm
was chosen much worse for the sake of robustness and reduction of computational complexity
(according to the concept of approximate computing) compared
to Gaussian method. Such a choice made the Richardson algorithm much faster than the Gaussian method.
In other words, optimization possibility of the Richardson algorithms associated with the trade-off
between accuracy on the one side  and robustness and computational time on the other side
allows the proper choice of acceptable uniform accuracy without any deterioration and essentially reduces the computational time.
\\
Finally, the Figure~\ref{figsurface3} shows that the approximation performance of the  Richardson algorithm
 is far superior than direct methods (which are not suitable for the
 detection of the power quality events) in the  ill-conditioning case.

\section{Recursive Parameter Calculation Algorithm}
\noindent
{\it Description.}
\noindent
The parameter vector in  (\ref{lineq1}) can be calculated using the
inverse of information matrix, $ \theta_k = A^{-1}_k  b_k$.
Denoting $ \Gamma_{k} = A^{-1}_k $ the recursive update of $ \Gamma_{k} $ via $ \Gamma_{k-1} $
is derived by application of the matrix inversion lemma\footnote{
$ (X + Y W Z)^{-1} = X^{-1} - X^{-1} Y ~ [ W^{-1} + Z~X^{-1}~Y]^{-1} ~ Z X^{-1} $ }, \cite{sto2023}
to the identity (\ref{ak1}):
\beq
 \Gamma_{k} = \Gamma_{k-1} - U_k ~ S^{-1} ~ U^T_k \label{gk1}
\eeq 	
where $Q_k = [\fai_{k} ~ \fai_{k-w}]$, $U_k =  \Gamma_{k-1} ~ Q_k$,
$S = D + Q^T_k ~ \Gamma_{k-1} ~ Q_k $ and  $ D = \begin{bmatrix}
1 & 0 \\
0 & -1  \end{bmatrix}.$ The $2 \times 2$ matrix $S$ remain the same in all the steps of the window of a given size $w$
and should be calculated only once.
Two forms of the parameter update $ \theta_k $ can be presented as follows:
\baq
\theta_k &=& [I - U_k ~ S^{-1} ~ Q^T_k] [ \theta_{k-1} + \Gamma_{k-1} d _k ] \label{tet11} \\
\theta_k &=& \Gamma_{k}   b_k \label{tet12}
\eaq
where $I$ is the identity matrix and the form (\ref{tet11}) is derived from (\ref{tet12}) and (\ref{gk1}).
The algorithms are initialized as  follows $\Gamma_{w} = A^{-1}_w$ and $A_w~\theta_w = b_w$.
The parameter update (\ref{tet12}) does not depend on parameters of the previous step and requires
matrix vector multiplication only. The inverse matrix
and the parameter update law,  (\ref{gk1}) and  (\ref{tet11}) can be calculated in two parallel loops.
Both forms are quadratic complexity algorithms and faster than direct parameter calculation  methods.
\\ {\it Error Accumulation.} The algorithm described above can be seen as ideal explicit recursive solution
of the system (\ref{lineq1}) - (\ref{fai1}).  Unfortunately, such solution is not robust with respect to error accumulation without corrections. The accumulation strength depends on the size of the moving window  and the information matrix.
Although the error accumulation is not very significant for relatively large window sizes the deterioration of the performance maybe essential for big data applications.
The performance deterioration due to error accumulation is significant for ill-conditioned information matrices and
for a large number of harmonics (which is expected in future electric networks) and short window sizes even for well conditioned information matrices due to the large number of calculations.
For the sake of improvement of the accuracy and robustness the Richardson corrections should be introduced in
(\ref{gk1}) and  (\ref{tet11}). \\
Notice that Richardson framework with  Newton-Schulz matrix inversion algorithms
is ideally suited for these corrections providing (after few iterations only) two improved estimates for the next step
of the algorithm.
\\ {\it Changeable Window Size.} Unfortunately, the algorithm (\ref{gk1}) and  (\ref{tet11}) should be re-initialized when the moving window changes the size, which happens quite often for the detection of both rapidly and slowly varying parameters. Initialization  includes calculations of the matrix inverse $\Gamma_{w} = A^{-1}_w$ and the parameter vector which satisfies  $A_w~\theta_w = b_w$ and it is computationally heavy for large scale systems.
Therefore new computationally efficient preconditioning methods should be developed for the case of frequent
changes of the window size.

\section{Properties of the Moving Window}
\label{proper}
\noindent
{\it Lemma.} Consider rank two update $R_k$ of the matrix $A_{k-1}$ defined in (\ref{ak1}).
Then eigenvalues of the matrix $A_k$ are the same for all the steps $k \ge w + 1$.
\\ {\it Proof.} The rank two matrix $R_k$ has two nonzero eigenvalues only,
$ \ds \pm \| R_k \|_F / \sqrt{2} $, where the norm is the Frobenius norm. To prove that the eigenvalues
remain the same it is sufficient to consider evolution of the coefficients of characteristic polynomial of this matrix
(or eigenvectors), \cite{sto2023}.
\begin{figure}
\centerline{\psfig{figure=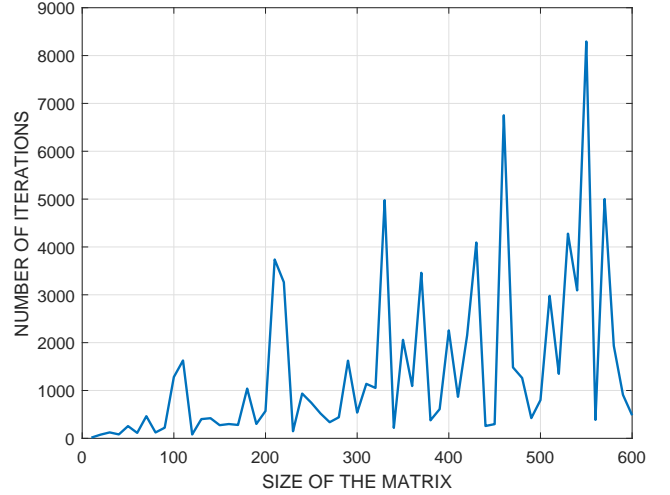,height=70mm}}
\begin{center}
\caption{\small{
The number of iterations (for estimation of the largest eigenpair $A x = \lambda x$)
as a function of the size of ill-conditioned information matrices
that is required to reach the following accuracy of estimation error:
$ \| \hat{x}_k ~ \|  A ~ \hat{x}_{k} \|  - A ~\hat{x}_{k} \| < 0.01 $. }}
\label{itern1}
\end{center}
\end{figure}
\section{Preconditioning Based on the Properties of the Window }
\noindent
{\it Simplest Preconditioner.}
The simplest preconditioner $\alpha = \frac{\ds 2}{\ds \| A \|_\infty} $  guarantees that
the spectral radius $\rho$  of the  SPD matrix $A$ is less than one, $\rho(I - \alpha A) < 1 $,
where $ \| \cdot \|_\infty $ is the maximum row sum matrix norm, \cite{ben1}, \cite{cdc54}.
\\
{\it Optimal Preconditioner \& Recursive Estimation of the Eigenvalues.}
The spectral radius of the matrix $(I - \alpha A)$ gets its minimal value
$(1 - \lambda_1 \alpha)$  for SPD matrix $A$ for the following preconditioner $ \alpha = \frac{\ds 2 }{ \ds \lambda_1 + \lambda_n} $, where $\lambda_1=\lambda_{min}(A) $ and $\lambda_n =\lambda_{max}(A)$ are minimal and maximal eigenvalues respectively. In other words the optimal preconditioner maps the interval which contains all eigenvalues of $A$ onto symmetric interval of maximal length around the origin, \cite{ben1}.
\\
The following power iteration algorithm \cite{din1}, (\ref{powm1}) which requires matrix vector multiplications only
\beq
 \hat{x}_k  = \frac{A ~ \hat{x}_{k-1}}{  \|  A ~ \hat{x}_{k-1} \| } \label{powm1}
\eeq
can be applied for estimation of the largest eigenpair $A x = \lambda x$.
\\ Notice that the minimal eigenvalue of $A$ can be estimated via maximal eigenvalue of $(\beta I - A)$
where $\beta = \hat{\lambda}_n + \epsilon$, $\hat{\lambda}_n$ is estimated maximal eigenvalue of $A$
and $\epsilon$  is a sufficiently small positive number.
Maximal eigenvalue of $(\beta I - A)$ in turn can be estimated using the same algorithm (\ref{powm1}).
\\
{\it Comparisons \& Drawbacks.} The spectral radius of the matrix $(I - \alpha A)$ with optimal preconditioner decreases only slightly  compared to the simplest one. In addition, numerical stability problems may occur in the algorithm
(\ref{powm1}) for estimation of the maximal eigenvalue in the presence of roundoff errors.
Estimation algorithms may require a large number of matrix vector multiplications for accurate estimation of
maximal eigenvalue and accurate estimation of minimal one may require even more matrix vector multiplications
due to error propagation. The number of iterations as a function of the size of the ill-conditioned information matrices that is required to reach the following accuracy of estimation error:
$ \| \hat{x}_k ~ \|  A ~ \hat{x}_{k} \|  - A ~\hat{x}_{k} \| < 0.01 $ is plotted in Figure~\ref{itern1}.
The Figure~\ref{itern1} shows that the number of iterations increases with the size of the matrix and can be sufficiently large (which corresponds to large number of matrix vector multiplications)
for large scale systems.
\\
 In addition, the spectral radius of the matrix $(I - \alpha A)$ with optimal preconditioning
can be even larger than one (which results in unstable system) for insufficient number of iterations of the
algorithm  (\ref{powm1}) (due to  insufficient computational capacity, for example) since the algorithm underestimates
the eigenvalue (due to convergence from below). On the contrary the simplest preconditioner which is
based on the maximum row sum matrix norm (associated with  the upper bound for all Gershgorin circles)
provides loose upper bound for maximal eigenvalue which is robust
in finite digit calculations and does not require any significant computational efforts.
\\ Notice that the upper bound of the largest eigenvalue can be estimated in many other ways, see
for example  \cite{haus}, \cite{fad} and  \cite{wol}. The choice of estimation algorithm
is associated with the trade-off between accuracy and computational complexity.
\\
{\it Suboptimal Preconditioner for Ill-Conditioned Cases.}
Minimal eigenvalue is very small for the ill-conditioned matrices and its accurate estimation requires significant computational efforts or even impossible for extreme ill-conditioning in finite digit calculations.
Small minimal eigenvalue can be neglected
in this case resulting in suboptimal preconditioner $ \alpha = \frac{\ds 2 }{ \ds \eps + \lambda_n} $ where
$\eps > 0$ is a sufficiently small number.
However, the spectral radius of the matrix $(I - \alpha A)$  with suboptimal preconditioner
maybe even larger then the spectral radius of the same matrix with the simplest preconditioner.
Indeed, the spectral radius of the matrix with the simplest preconditioner is associated with
minimal eigenvalue $\lambda_1$  and is calculated as  $1 -  \frac{\ds 2}{\ds \| A \|_\infty} \lambda_1$,
whereas the spectral radius of suboptimal preconditioner is the absolute value of
$1 - \frac{\ds 2 \lambda_n}{\ds \lambda_n + \eps} $ which depends on the maximal eigenvalue $\lambda_n$ only
and can be closer to one.
\\ {\it Recommendations.} Optimal and suboptimal preconditioners can be applied for the case where
sufficient computational capacity is available in preprocessing and the window size does not change
during processing. Then the preconditioner that is based on estimated largest eigenvalue can be applied in all the steps since the eigenvalues  are the same, see Lemma in Section~\ref{proper}
\footnote{Notice that estimation of the largest eigenvalue associated with changeable window size can also be performed using parallel computational units (or in memory) and send to the signal processing unit}. Otherwise the simplest
preconditioner which does not require any computational  efforts (compared to optimal preconditioner)
can be applied.

\section{Conclusions}
\label{conc}
 \noindent
Application of the direct methods for parameter calculations (as solution of the algebraic equations)
provide inaccurate results in ill- conditioned case due to divisions by very small numbers, which has direct negative impact on estimation performance as it is shown on real data in this report.
Richardson approach does not have any divisions by small numbers in the ill-conditioned case
and provides robust and accurate solution.
\\ It was also shown that simple preconditioner, which is based on maximum row sum matrix norm
provides more robust results than optimal preconditioner based on estimation of the eigenvalues.




\end{document}